\documentclass[11pt, reqno]{amsart}
\usepackage{amsmath, amsthm, amscd, amsfonts, amssymb, graphicx, color}
\usepackage[bookmarksnumbered, colorlinks, plainpages]{hyperref}

\usepackage[left=2.5 cm, right= 2.5 cm, top = 3 cm, bottom = 2.5 cm]{geometry}

\def\R{\mathbb{R}}

\def\Z{\mathbb{Z}}
\def\N{\mathbb{N}}

\def\lt{\left}
\def\rt{\right}

\newtheorem{theorem}{Theorem}

\theoremstyle{definition}
\newtheorem{definition}[theorem]{Definition}

\theoremstyle{remark}
\newtheorem{remark}[theorem]{Remark}
\begin{document}
\setcounter{page}{1}

\centerline{}

\centerline{}

\vspace{-20pt}

\title[A. Prade]{Boundary H\"older regularity for the fractional Laplacian \\ over Reifenberg flat domains via ABP maximum principle}

\author[A. Prade]{Adriano Prade}

\address{CMAP, \'Ecole Polytechnique, 91120 Palaiseau, France.}
\email{adriano.prade@polytechnique.edu}

\centerline{ \large Springer INdAM Series, Conference ``Geometric Measure Theory and applications 2024''}

\vspace{15 pt}

\begin{abstract} For $0<s<1$, we consider the nonlocal equation $(-\Delta)^s u = f$ over a Reifenberg flat domain $\Omega$ with $f \in C({\overline{\Omega}})$ and null Dirichlet exterior condition. Given $\alpha \in (0,s)$, we prove that weak solutions are $\alpha$-H\"older continuous up to the boundary when the flatness parameter is small enough. The main ingredients of the proof are an iterative argument and a nonlocal version of the ABP maximum principle. 
\end{abstract} 

\maketitle

The problem of boundary regularity of elliptic partial differential equations over Reifenberg flat sets has gained some attention in recent years, both in the local and in the nonlocal case. In \cite{lemenant2013boundary}, Lemenant and Sire proved boundary H\"older regularity for the Poisson equation with $L^q$ datum and null Dirichlet boundary condition. %, by using a variant of Alt-Caffarelli-Friedman monotonicity formula and Campanato criterion.
Later on, Lian and Zhang extended the previous result to more general classes of elliptic equations, see \cite{lian2018boundary}. They considered linear equations in divergence form, the $p$-Laplace equation and even fully nonlinear equations in non-divergence form. %The main tools behind their proof are an induction argument exploiting the structure of Reifenberg flat sets and the classical ABP maximum principle (see \cite{caffarelli1995fully} for the general fully nonlinear version or \cite{gilbarg1977elliptic} for the linear case).

More in general, the regularity of nonlocal elliptic equations has been a popular research topic in the last decade and many major contributions are reported in the the monograph \cite{fernandez2024integro} by Fern\'andez-Real and Ros-Oton. For example, boundary H\"older regularity was proven for equations over $C^{1,\gamma}$, $C^1$ and Lipschitz domains. Afterwards, Goldman, Novaga and Ruffini showed in \cite[Lemma 3.17]{goldman2024rigidity} a boundary H\"older estimate for the potential function of a Reifenberg flat set $\Omega \subset \R^n$. %, by adapting the arguments of Lemenant and Sire.
Since the function they considered is the solution to a nonlocal equation involving $(- \Delta)^{\frac{1}{2}}$, their conclusion marked a first step towards boundary H\"older regularity for the fractional Laplacian on Reifenberg flat sets. \\

The purpose of this proceeding is to complete the picture for any $0<s<1$, by proving the same result for general distributional solutions to the system:
\[ \begin{cases}
        (-\Delta)^s u = f & \mbox{on } \Omega \\ 
        u=0 & \mbox{on } \mathbb{R}^n \setminus \Omega,
    \end{cases} \]
where $\Omega$ is a $(\eta, r_0)$-Reifenberg flat set with parameter of flatness $\eta$ small enough and $f \in C({\overline{\Omega}})$. The proceeding is organized as follows. We recall the definition of Reifenberg flatness, then we state our boundary H\"older regularity result and we outline the general strategy of the proof. We introduce a suitable barrier for the problem and finally we prove the theorem. Hence, we begin by presenting the definition of a $(\eta, r_0)$-Reifenberg flat set.
\begin{definition}
Let $ 0< \eta < 1/2 $ and $r_0 > 0$. A set $\Omega \subset \mathbb{R}^n$ is $(\eta, r_0)$-Reifenberg flat if for all $x \in \partial \Omega$ and for all $0<r<r_0$ there exists a hyperplane $H_{x,r}$ containing $x$ and such that:
    \begin{itemize}
        \item $d_H ( \partial \Omega \cap B_r(x), H_{x,r} \cap B_r(x) ) \leq \eta r$, where $d_H$ denotes the Hausdorff distance;
        \item one of the connected components of $\{d(\cdot, H_{x,r}) \geq 2\eta r \} \cap B_{r}(x)$ is included in $\Omega$ and the other in $\Omega^c$.
    \end{itemize}
\end{definition}
Essentially, we can see the boundary of a Reifenberg flat set as a perturbation of a hyperplane at sufficiently small scales. \\

Now we state our boundary H\"older regularity result for the fractional Laplacian. In the following, we always refer to \cite{fernandez2024integro} for any preliminary result regarding the fractional Laplacian and nonlocal equations, for example the different notions of solution, their properties and the relations between them. Also, from now on we assume $n\geq 2$ and we denote $d_{\Omega}(x) := \text{dist}(x, \Omega^c)$ for any $x\in \Omega$.
\begin{theorem}\label{th}
   Let $s \in (0,1)$ and $\alpha \in (0,s)$. Let $\Omega \subset \R^n$ be a bounded, open and connected $(\eta, r_0)$-Reifenberg flat domain, $f \in C(\overline{\Omega})$ and $u$ be the weak solution to the system
   \begin{equation}\label{system} \begin{cases}
        (-\Delta)^s u = f & \text{on } \Omega \\ 
        u=0 & \text{on } \mathbb{R}^n \setminus \Omega.
    \end{cases} \end{equation}
   Then, there exists a constant $\eta_0(n,s,\alpha) >0$ such that for all $\eta \leq \eta_0$ we have $u \in C^{\alpha}(\overline{\Omega})$, with the estimate
   \begin{equation}\label{holdercontinuity}
       \|u\|_{C^{\alpha}(\overline{\Omega})} \leq C \| f\|_{L^{\infty}(\Omega)}
   \end{equation}
   for some constant $C=C(n,s,\alpha)>0$.
\end{theorem}
%Existence of a unique weak solution $u \in L^{\infty}(\R^n)$ for the above problem is ensured by \cite[Theorem 2.2.24 and Lemma 2.3.9]{fernandez2024integro}, it is also a distributional solution by \cite[Lemma 2.2.32]{fernandez2024integro}
We already know that $u \in L^{\infty}(\R^n)$ thanks to the $L^{\infty}$ bound for weak solutions \cite[Lemma 2.3.9]{fernandez2024integro}. In addition, $u$ is also a distributional solution to \eqref{system} by the equivalence result \cite[Lemma 2.2.32]{fernandez2024integro}. The proof of Theorem \ref{th} follows the general strategy employed for $C^{1,\gamma}$, $C^1$ and Lipschitz domains in \cite[Section 2.6]{fernandez2024integro}. We first show that
\begin{equation}\label{estimatedistance}
    |u(x)| \leq C d_{\Omega}^{\alpha}(x) \qquad \text{for } x \in \Omega \text{ and } C>0. 
\end{equation}
Then, we combine it with interior regularity estimates and we conclude that $u \in C^{\alpha}(\overline{\Omega})$. Estimates like \eqref{estimatedistance} are usually obtained via comparison principle by using a barrier function with the proper $\alpha$-H\"older growth. The construction of a barrier depends on the regularity of the domain $\Omega$ where the equation is defined and some of them are suitable even for Reifenberg flat sets, see \cite[Lemma B.3.3]{fernandez2024integro}. However, in this case we have too little information about the boundary, so it is not easy to follow other arguments valid for more regular domains (see for example \cite[Proposition 2.6.4]{fernandez2024integro} for the $C^{1, \gamma}$ case).

A similar but more convenient approach was adopted in \cite{lian2018boundary} to prove boundary H\"older regularity for many classes of local elliptic equations. Since Reifenberg flat sets can be approximated by hyperplanes at every scale around every point of their boundary, the key idea is to use an iterative argument to get progressively more information on the behaviour of the solutions. Practically speaking, the authors first introduced a barrier with respect to a flat boundary and then proceeded by induction, rescaling it and applying the ABP maximum principle. In this way, they obtained the desired estimate on the solutions at every step. \\

To prove Theorem \ref{th}, we adapt this method to our nonlocal situation. With the help of some potential theory, we first build a proper barrier function $v_0$. Then, after rescaling, we apply Guillen and Schwab's nonlocal ABP maximum principle from \cite{guillen2012aleksandrov} to draw the requested comparison. \\

Denoting $\{e_i\}_{i=1}^n$ the standard basis of $\R^n$, we consider the closed ball $\Tilde{B_1} := \overline{B \lt(- \lt(\frac{1}{8} +\eta \rt) e_n, \frac{1}{8} \rt)}$. The $2s$-Riesz energy $\mathcal{I}_{2s}(\Tilde{B}_1)$ of $\Tilde{B_1}$ is defined as:
\begin{equation}\label{riesz}
    \mathcal{I}_{2s}(\Tilde{B}_1) := \min_{\mu(\Tilde{B}_1)=1} \int_{\Tilde{B}_1 \times \Tilde{B}_1} \frac{d\mu(x)d\mu(y)}{|x-y|^{n-2s}}.
\end{equation}
There exists an unique probability measure $\mu_{\Tilde{B}_1}$ attaining the minimum in \eqref{riesz}, see \cite{landkof1972foundations}. We define the potential function $v_{\Tilde{B}_1} \! \!: \R^n \longrightarrow \R $ as:
\[ v_{\Tilde{B}_1}(x) := \int_{\Tilde{B}_1} \frac{d\mu_{\Tilde{B}_1}(y)}{|x-y|^{n-2s}}. \]
We take as a barrier for our proof the function $v_0 : \R^n \longrightarrow \R $ defined as:
\[ v_0(x) := 1 - \mathcal{I}_{2s}(\Tilde{B}_1)^{-1} v_{\Tilde{B}_1}(x). \]
From the properties of $v_{\Tilde{B}_1}$ presented in \cite[Section 2]{goldman2015existence}, we know that $v_0 \in [0,1]$, it is radially symmetric and monotone increasing around the point $-\lt(\frac{1}{8}+\eta\rt)e_n$. In addition, the system
    \begin{equation}\label{barrier}
        \begin{cases}
            (-\Delta)^s v_0 = 0 & \text{on } \R^n \setminus \Tilde{B}_1 \\ v_0 = 0 & \text{on } \Tilde{B}_1 
            % \\  v_0(x) \longrightarrow 1 & \text{as } |x| \rightarrow{ + \infty}
        \end{cases}
    \end{equation}
holds pointwise, so also in the distributional sense. It is possible to show that $v_0 \in C^s(\R^n)$ and we denote $C_H=C_H(n,s)$ the $s$-H\"older constant. In particular, for all $\eta >0$ we get the estimate:
    \begin{equation}\label{holderbarrier}
        v_0(x) \leq C_H\eta^s \qquad \text{when} \quad d(x, \tilde{B}_1) \leq 2 \eta.
    \end{equation}
We also have $v_0(x) \geq d(x, \tilde{B}_1)^s$, thanks to the maximum principle for the fractional Laplacian \cite[Lemma 1.10.8]{fernandez2024integro}. Finally, with our choice of the set $\tilde{B}_1$ it is easy to show $d(x, \tilde{B}_1) \geq \frac{1}{4}|x|$ whenever $|x|>1$, therefore:
    \begin{equation}\label{major}
        v_0(x) \geq \frac{1}{4^s}|x|^s \qquad \text{for} \quad |x|>1.
    \end{equation}
We are now ready to prove Theorem \ref{th}.    
\begin{proof}
    The statement is invariant under dilation and translation, so we assume $r_0=1$ and $x=0$ without loss of generality. We set $M= \| u \|_{L^{\infty}(\Omega)} + C_{\text{ABP} } \,\| f\|_{L^{\infty}(\Omega)}$, where $C_{\text{ABP}}$ appears later in \eqref{firstABP}. Our goal is to show that there exist $0 < \eta (n, s, \alpha) < 1/2$ and $0 < \lambda (\eta) <1$ such that:
    \begin{equation}\label{induction}
        \| u \|_{L^{\infty}(\Omega \, \cap B_{\lambda^k})} \leq C M \lambda^{k \alpha} \qquad \text{for all } k \in \Z,
    \end{equation}
    for another constant $C=C(n,s,\alpha) \geq 1$ which is explicitly determined afterwards. We prove \eqref{induction} by induction. First, if $k \in \Z$ and $k \leq 0$ the statement holds, since $ u \in L^{\infty} (\R^n)$, $C\geq1$ and $0 < \lambda <1$. Then, we fix $k\in \N$ and we assume that:
    \begin{equation}\label{inductionhypo}
    \| u \|_{L^{\infty}(\Omega \, \cap B_{\lambda^j})} \leq C M \lambda^{j \alpha} \qquad \text{for all } j \in \Z \text{ with } j \leq k.
    \end{equation}
    We show that there exists a number $m \in \N$ (depending only on $s$ and $\alpha$) such that:
    \begin{equation}\label{conclusion}
    \| u \|_{L^{\infty}(\Omega \, \cap B_{\lambda^{k+m}})} \leq C M \lambda^{(k+m)\alpha}. \end{equation}
    To do so, we work on the set $\Omega \, \cap B_{\lambda^k}$ and take as barrier $v_k$ a suitable rescaling of $v_0$:
    \[ v_k(x) = 4^s CM\lambda^{(k-1)\alpha} v_0 \lt( \frac{x}{\lambda^k} \rt). \]
    Denoting $\Tilde{B}_{\lambda^k} := \lambda^k \Tilde{B_1}$ we automatically get from \eqref{barrier}: %$v_k \in [0,CM(\lambda^{k-1})^{\alpha}]$:
    \begin{equation}\label{rescaledbarrier}
        \begin{cases}
            (-\Delta)^s v_k = 0 & \text{on } \R^n \setminus \Tilde{B}_{\lambda^k} \\ v_k = 0 & \text{on } \Tilde{B}_{\lambda^k}. 
            %\\ v_k(x) \longrightarrow 4^s CM(\lambda^{k-1})^{\alpha} & \text{as } |x| \rightarrow{ + \infty}
    \end{cases} \end{equation}
    In order to apply Guillen and Schwab's ABP maximum principle we have to show that the system
    \begin{equation}\label{viscositysystem}
        \begin{cases}
            (-\Delta)^s (v_k - u) = - f & \text{on } \Omega \, \cap B_{\lambda^k} \\ (v_k - u) \geq 0 & \text{on } \lt(\Omega \, \cap B_{\lambda^k} \rt)^c
    \end{cases}
    \end{equation}
    holds in viscosity sense, see \cite[Section 3.2]{fernandez2024integro} for more details about viscosity solutions to nonlocal equations. The first condition holds in distributional sense and it is an easy consequence of \eqref{system} and \eqref{rescaledbarrier} once noticed that $ \Omega \cap B_{\lambda^k} \subset \R^n \setminus \Tilde{B}_{\lambda^k}$. By the equivalence between distributional and viscosity solutions in the linear case \cite[Lemma 3.4.13]{fernandez2024integro}, the first line of \eqref{viscositysystem} is true in viscosity sense too since $f \in C(\overline{\Omega})$. Now we check that $v_k - u\geq 0$ over the set $\lt(\Omega \cap B_{\lambda^k} \rt)^c$: we immediately get it for $x \in \Omega^c$ because $u(x)=0$ and $v_k \geq 0$ on $\R^n$. Conversely, if $x \in \Omega \setminus B_{\lambda^k}$ then
    \[
    |x| \in [ \lambda^j, \lambda^{j-1}) \qquad \text{for some   } j\in \Z \text{ with } j \leq k.
    \]
    By induction hypothesis \eqref{inductionhypo}, we have $ u(x) \leq C M \lambda^{(j-1)\alpha}$. Moreover, since $|x|> \lambda^k$ we use \eqref{major} and we also get:
    \[ v_k(x) = 4^s CM\lambda^{(k-1)\alpha} v_0 \lt( \frac{x}{\lambda^k} \rt) \geq CM\lambda^{(k-1)\alpha}  \frac{|x|^s}{\lambda^{ks}}  \geq CM\lambda^{(k-1)\alpha} \lambda^{(j-k)s}. \]
    Altogether:
    \[ \frac{v_k(x)}{u(x)} \geq \lambda^{(j-k)(s- \alpha)} \geq 1 \qquad \text{since} \quad j-k \leq 0,\quad s-\alpha >0\, \quad\text{and}\quad 0<\lambda <1;\]
    so we finally obtain $v_k - u\geq 0$ on the whole of $\lt(\Omega \cap B_{\lambda^k} \rt)^c$ as well as the validity of \eqref{viscositysystem}. Therefore, by Guillen and Schwab's ABP maximum principle we get:
    \begin{equation}\label{firstABP}
        \sup_{\quad\Omega \cap B_{\lambda^k}} (v_k - u)^- \leq C_{\text{ABP}} \,\,\text{diam}(\Omega \cap B_{\lambda^k}) \,\| f\|_{L^{\infty}(\Omega \cap B_{\lambda^k})} \leq M \lambda^k.
    \end{equation}
    We may rewrite it as
    \[ u(x) \leq 4^s CM\lambda^{(k-1)\alpha} v_0 \lt( \frac{x}{\lambda^k} \rt) + M\lambda^k \qquad \text{for all} \,\, x \in \Omega \cap B_{\lambda^k}. \]
    \\
    For all $x \in \Omega \cap B_{\eta \lambda^k}$ there holds $d(x, \tilde{B}_{\lambda^k}) \leq 2 \eta \!\cdot\! \lambda^k$, so we have by \eqref{holderbarrier}:
    \begin{equation}\label{ABP}  
    u(x) \leq 4^s CMC_H\lambda^{(k-1)\alpha} \eta^s + M\lambda^k \qquad \text{for all} \,\, x \in \Omega \cap B_{\eta \lambda^k}.
    \end{equation}
    We are ready to select the right parameters $m$, $\lambda$ and $C$ to get \eqref{conclusion}. Since $s - \alpha > 0$, there exists a number $m=m(s,\alpha) \in \N$ such that $m(s-\alpha) > s$, or $\frac{s + m\alpha}{s} < m$ equivalently. We set 
    \[ \eta = \lambda^{\frac{s + m\alpha}{s}}. \]
    As a consequence $\lambda^m < \eta$, so we get $B_{\lambda^{k+m}} \subset B_{\eta \lambda^k}$. Hence, from \eqref{ABP} we obtain for all $x \in \Omega \cap B_{\lambda^{k+m}}$:
    \[ \begin{split}
        u(x) & \leq 4^s CMC_H\lambda^{(k-1)\alpha} \lambda^{s+m\alpha} + M\lambda^k \\ & = CM \lambda^{(k+m)\alpha} \lt( 4^s C_H \lambda^{s-\alpha} + \frac{1}{C} \lambda^{k-(k+m)\alpha} \rt).
    \end{split}\]
    Now, we choose first $\lambda$ small enough such that there hold at the same time $4^s C_H \lambda^{s-\alpha} \leq 1/2$ and $\eta(\lambda) < 1/2$. Next, we take $C$ large enough such that $C \geq 1$ and $\frac{1}{C}\lambda^{k-(k+m)\alpha} \leq 1/2$. Notice that $k - (k+m)\alpha \rightarrow + \infty$ as $k \rightarrow + \infty$, so the choice of the constant $C$ can be made independently of $k \in \N$. In this way, we get:
    \[
    u(x) \leq CM \lambda^{(k+m)\alpha} \qquad \text{for all} \,\, x \in \Omega \cap B_{ \lambda^{k+m}}.
    \]
    To complete the proof of \eqref{conclusion} we only need to show the lower bound:
    \[ u(x) \geq - CM \lambda^{(k+m)\alpha} \qquad \text{for all} \,\, x \in \Omega \cap B_{ \lambda^{k+m}}. \]
    The argument is the same as before, with the only difference that in this case we apply the ABP maximum principle on the function $u + v_k$, so that in the end we get $u \geq -v_k - M \lambda^k$. Therefore, we obtain \eqref{conclusion} and we conclude \eqref{induction} by induction. In particular, for all $x\in \Omega$ and $y \in \partial \Omega$ we have:
    \[ |u(x)| \leq CM |x-y|^{\alpha}.\]
    Choosing $x\in B_{3/4}$ and $\overline{x} \in \partial \Omega$ such that $|x-\overline{x}|= d_{\Omega}(x)$, we also get the estimate:
    \begin{equation}\label{pointholder}
        |u(x)| \leq CM d_{\Omega}^{\alpha}(x) \qquad \text{for all} \,\, x \in B_{3/4}.
    \end{equation}
    By arguing as \cite[Proposition 2.6.4, Step 2]{fernandez2024integro}, we combine \eqref{pointholder} with interior regularity estimates and we deduce that $u \in C^{\alpha}(B_{1/2})$. Finally, we cover $\Omega$ with finitely many balls of radius $1/2$ and we reason in the same way on each of them thus obtaining $u \in C^{\alpha}(\overline{\Omega})$. Estimate \eqref{holdercontinuity} comes directly from \eqref{induction} and the $L^{\infty}$ bound for weak solutions of \cite[Lemma 2.3.9]{fernandez2024integro}.
\end{proof}
\begin{remark}
    In the proof we apply a rather unusual version of an induction argument, where the inductive step is obtained thanks to a natural number $m$ different than $1$. Indeed, assuming the induction hypothesis \eqref{inductionhypo}, we cannot directly conclude 
    \[ u(x) \leq CM \lambda^{(k+1)\alpha} \qquad \text{for all} \,\, x \in \Omega \cap B_{ \lambda^{k+1}}, \] 
    because we want to get a bound on $u$ which is also consistent with the set where it holds. To be more precise, we start from the estimate \eqref{ABP} and we point out the following observations.
    \begin{enumerate}
        \item One could be tempted to set $\eta = \lambda$ (as it is successfully done in \cite{lian2018boundary} for the classical Laplacian case), in order to immediately get an estimate on $\Omega \cap B_{ \lambda^{k+1}}$. However, in this way we would get:
        \[ u(x) \leq 4^s CM\lambda^{(k+1)\alpha} C_H \lambda^{s-2\alpha} + M\lambda^k \qquad \text{for all} \,\, x \in \Omega \cap B_{ \lambda^{k+1}}, \]
        which yields H\"older regularity only for the values of $\alpha$ for which $s-2\alpha > 0$, namely $\alpha \in (0, s/2)$.
        \item To avoid this issue, we could find the right value of $\eta$ such that there holds $\lambda^{(k-1)\alpha} \eta^s = \lambda^{(k+1)\alpha} \lambda^{s-\alpha}$, it corresponds to $\eta = \lambda^{\frac{s+\alpha}{s}}$. Anyway, this choice implies $\lambda^2 < \eta < \lambda$, so we can only conclude 
        \[ u(x) \leq 4^s CM\lambda^{(k+1)\alpha} C_H \lambda^{s-\alpha} + M\lambda^k \qquad \text{for all} \,\, x \in \Omega \cap B_{ \lambda^{k+2}}, \]
        but not for all $x \in \Omega \cap B_{ \lambda^{k+1}}$.
    \end{enumerate} 
\end{remark}

The boundary regularity result over Reifenberg flat domains we presented may be extended to more general nonlocal elliptic operators and it will be the subject of the upcoming work \cite{prade2025boundary}. \\

{\bf Acknowledgement.} We kindly thank M. Goldman for very helpful discussions leading to this result and general advice. 

\bibliographystyle{acm}
\bibliography{Proceeding}

%its $2s$-Riesz energy is defined as:
%\[ \mathcal{I}_{2s}(\Tilde{B}_1) := \min_{\mu(\Tilde{B}_1)=1} \int_{\Tilde{B}_1 \times \Tilde{B}_1} \frac{d\mu(x)d\mu(y)}{|x-y|^{n-2s}}\]
%Given the measure $\mu_{\Tilde{B}_1}$ attaining the minimum, it is well defined the potential function $v_{\Tilde{B}_1} \! \!: \R^n \longrightarrow \R $ as:
%\[ v_{\Tilde{B}_1}(x) := \int_{\Tilde{B}_1} \frac{d\mu_{\Tilde{B}_1}(y)}{|x-y|^{n-2s}} \]

\end{document}